\documentclass[12pt]{article}
\usepackage{amssymb,latexsym, amsmath, amscd, array, graphicx}
\usepackage{amssymb}

\usepackage{setspace}
\input xy
\xyoption{all}
\makeatletter 

\usepackage{setspace}
\usepackage{enumitem}
\date{}

\newtheorem{theorem}{Theorem}[section]

\newtheorem{proposition}[theorem]{Proposition}

\newtheorem{corollary}[theorem]{Corollary}

\newtheorem{problem}[theorem]{Problem}

\newcommand{\z}{{\Bbb Z}}

\newcommand{\re}{{\Bbb R}}

\newcommand{\invlim}{\raisebox{-1ex}{$\stackrel{\hbox{lim}}{\leftarrow}$}}

\newcommand{\lo}{\rightarrow}

\newcommand{\p}{\dim_{\z[\frac{1}{p}]} }

\newcommand{\black}{{\blacksquare}}

\begin{document}

\title{\bf On $p$-adic actions raising  dimension   by $2$}

\author{  Michael  Levin\footnote{This research was supported by 
THE ISRAEL SCIENCE FOUNDATION (grant No. 522/14) }}

\maketitle
\begin{abstract} Raymond and Wiliams 
constructed an action of  the $p$-adic integers $A_p$ on an $n$-dimensional compactum $Z$, $n \geq 2$, with 
the orbit space $Z/A_p$ of dimension $n+2$. A simplified  construction of such an action was presented in
\cite{levin-on-rw}. In this paper we generalize the approach of \cite{levin-on-rw} to show that
 for every $n \geq 2$, an $(n+2)$-dimensional compactum $X$ can be obtained as the orbit space
$X=Z/A_p$ of an action of $A_p$ on an $n$-dimensional compactum $Z$ if and only if 
$\p X \leq n$ where $\p X$ is the cohomological dimension of $X$ with  coefficients in
$\z[\frac{1}{p}]$.
\bigskip
\\
{\bf Math. Subj. Class.:}  55M10, 22C05 (54F45)

\end{abstract}
\begin{section}{Introduction} 
Interest in (continuous) actions of the $p$-adic integers $A_p$  on compacta (=compact metric spaces) is inspired by 
the Hilbert-Smith conjecture that asserts that a compact group acting effectively on a manifold  must be a Lie group.
This conjecture is equivalent to the following  one: the group  $A_p$  of the $p$-adic integers cannot act effectively 
on a manifold. Yang \cite{yang} showed  that if $A_p$ acts effectively on a manifold $M$ then
either $\dim M/A_p=\infty$ or 
$\dim M/A_p =\dim M +2$.
 In order to verify if the latter dimensional relation ever occurs in a more general setting
Raymond and Williams \cite{raymond-williams} constructed an action of $A_p$  on an $n$-dimensional 
compactum (=compact metric space) $Z, n \geq 2$, 
with $\dim Z/A_p=n+2$. The author presented in \cite{levin-on-rw}  a simpler construction of  such an example.  
In this paper we generalize the approach of \cite{levin-on-rw} to provide  a full characterization of 
finite-dimensional compacta $X$
with $\dim X \geq 4$ 
that can be obtained as the orbit spaces $X=Z/A_p$ of an action of $A_p$ on a compactum $Z$ with $\dim X=\dim Z+2$.

Recall that the cohomological dimension of a compactum $X$ with respect to an abelian coefficient group $G$
is the least integer $n$ such that  the Cech cohomology $H^{n+1}(X, F; G)$ vanishes  for every closed
$F  \subset X$.
Clearly $\dim_G X \leq \dim X$ if $X$ is finite-dimensional. It is well-known that
for an action of $A_p$ on a compactum $Z$ we have that $\p Z =\p Z/A_p$.  Thus every 
$(n+2)$-compactum $X$ that can be represented as the orbit space of an action of $A_p$ on a $n$-dimensional
compactum  must satisfy $\p X  \leq n$. The goal of this paper is to show that the  inequality 
$\p X \leq  n$ is not only
necessary but also sufficient for   such a representation of $X$.

\begin{theorem}
\label{main-th}
Let $X$ be an $(n+2)$-dimensional  compactum with $n\geq 2$. Then $X$ can be obtained as the orbit space
$X=Z/A_p$  of
an action of $A_p$ on an $n$-dimensional compactum $Z$  if and only if $\p X \leq n$.
\end{theorem}

The paper is organized as follows. Sections 1 and  2 are devoted to constructions and properties 
needed  for proving Theorem \ref{main-th} in Section 3.
A few related remarks are given in Section 4.

\end{section}
\begin{section}{Auxiliary constructions and properties}
Consider  the unit $(m+2)$-sphere  $S^{m+2}$ in $\re^{m+3}, n\geq 2$. Represent $\re^{m+3}$ as 
the product $\re^{m+1} \times \re^2_\perp$
 of a coordinate plane $\re^{m+1}$ and its orthogonal complement $\re^2_\perp$ and let $S^m=S^{m+2} \cap \re^{m+1}$ 
 and
 $S^1_\perp =S^{m+2} \cap \re^2_\perp$. Take the  closed $\epsilon$-neighborhood $F$ of $S^m$.
 Then $F$ and  $F_\perp=S^{m+2} \setminus {\rm int} F$ can be represented as the products
 $F=S^m\times D$  and 
  $F_\perp= D_\perp \times  S^1_\perp $ of an $m$-sphere  $S^m$ and a circle  $S^1_\perp$ with
 a $2$-disk  $D$  and  an $(m+1)$-ball $D_\perp$  respectively. Thus we have that 
 $\partial F = \partial F_\perp =S^m \times \partial D = \partial D_\perp \times S^1_\perp$.
  Let $T=\re  /\z$ freely act on the circle $\partial D$ by rotations. Then this action induces  
the corresponding free action  on  $\partial F  = \partial F_\perp$ and the later action obviously extends  over 
$  F_\perp  $ as a free action and  it extends over $F$  by the  rotations of  $D$ induced by the rotations of $\partial D$.
This way we have defined an action of $T$ on $S^{m+2}$ whose  fixed point set  is $S^m$ and
$T$ acts freely on $S^{m+2}\setminus S^m$. Note that for every subgroup 
$\z_{p^k}=\z/p^k\z$ of $T$ the space $S^{m+2}/\z_{p^k}$ is homeomorphic to an
$(m+2)$-sphere.

 \begin{proposition}
\label{prop0}
Let $\Delta$ be an $(m+3)$-simplex. There is an unknotted sphere  $S^{m}$ in $S^{m+2}=\partial \Delta$ such that
for each $(m+2)$-simplex $\Delta'$ the intersection $S^m \cap \Delta'$ is an  $m$-ball whose  boundary  is  in 
$\partial \Delta'$ and whose  interior is in the interior of $\Delta'$.
\end{proposition}
{\bf Proof.}
The case $m=1$ was proved in \cite{levin-on-rw}.  We will proceed by induction on $m$.
Fix an $(m+2)$-simplex $\Delta_*$ of $\Delta$ and let $v$ be the vertex of $\Delta$ that does not belong to $\Delta_*$
and let $v_*$ be the barycenter of $\Delta_*$. Consider a sphere $S^{m-1}$ embedded in $\partial \Delta_*$
as required in the proposition with $m$ replaced by $m-1$ and $\Delta$ by by $\Delta_*$.
Then the union of the cones over the sphere $S^{m-1} \subset \Delta$ with the vertices $v_*$ and $v$ respectively  forms
a  sphere $S^m$ that satisfies the conclusions of the proposition.
 $\black$
\\
\\
From now we identify $S^{m+2}$ with the boundary of $(m+3)$-simplex $\Delta$ as described in
Proposition \ref{prop0}.
 Recall that 
 $F$  is represented as  the product $F=S^m\times D$, take $(m+4)$ disjoint closed $m$-balls
 $B'_1$, ..., $B'_{m+4}$ in $S^{m}$ and  denote  $B_i =B'_i \times D$. Then each closed $(m+2)$-ball 
 $B_i$ is invariant under the action of $T$
 and hence $M=S^{m+2}\setminus {\rm int}(B_1\cup\dots\cup B_{m+4})$ is invariant as well. 
 By Proposition \ref{prop0} we may assume that
  the balls $B_i$  lie in the interior of different $(m+2)$-simplexes $\Delta_i$ of $\Delta$.

\begin{proposition}
\label{prop1}
${}$

(i) For every $g \in T$ the map $x \lo gx, x \in M,$  is ambiently   isotopic
to a homeomorphism $\phi : M \lo M$ 
 by an isotopy of $M$  that does not move the points of $S^m\cup \partial  M$ 
and such that  $\phi$
restricted to the $m$-skeleton $\Delta^{(m+1)}$ of $\Delta$ coincides with
 the inclusion of $\Delta^{(m+1)}$ into $M$ and $\phi(\Delta' \cap M)=\Delta' \cap M$ for
 every $(m+2)$-simplex $\Delta'$ of $\Delta$.

(ii)  There is a retraction of $r: M\lo \Delta^{(m+1)}$ such that for every $(m+2)$-simplex $\Delta'$ of $\Delta$
we have  that $r(\Delta'  \cap M) \subset \partial \Delta'$.
\end{proposition}
{\bf Proof.} 

(i) Since $T$ is path-connected $g : M \lo M$ is isotopic to the identity map of $M$ by an isoptopy that does not
move the points of $M\cap S^m$ (the fixed point set of $T$). For each $\partial B_i \subset \partial M$ we can
change this isotopy on a neighborhood of $\partial B_i$ in $M$ to get in addition that the points of
$\partial B_i$ are not moved.

(ii) Define $r$ as on each $(m+2)$-simplex $\Delta'$ of $\Delta$  as  the restriction of a radial projection to
$\partial \Delta'$ from a point 
of $\Delta'$ that
that does not belong to $M$.
$\black$
\\ \\
For a prime number  $p$  and a positive integer $k$ 
consider the subgroup $\z_{p^k}$ of $T$ and  let 
 $\Gamma_M  =\Gamma_M(p^k) =M/\z_{p^k}$.   Note that $\Gamma_M$ is homeomorphic to $M$.
 Denote by $M_+$ the space which is the union of $p^k$ copies $M_a$, $a \in \z_{p^k}$, of $M$
 with $\partial M$ being identified in all the copies by the identity map. Thus
 all $M_a \subset M_+$ intersect each other at $\partial M$.  Define the action of $\z_{p^k}$ on $M_+$
 by sending  $x \in M_a$ to $gx \in M_{g+a}$.  Note that $\partial M \subset M_+$ is invariant under the action of
 $\z_{p^k}$ on $M_+$ and the natural  projection  $\alpha_+ : M_+ \lo M$ sending each $M_a$ to $M$ by the identity map
 is  equivariant.
 
 Let $\Gamma_+=\Gamma_+(p^k)$ 
 be the mapping cylinder of $\alpha_+$. The actions of $\z_{p^k}$ on $M_+$ and $M$
 induce the corresponding action of $\z_{p^k}$ on $\Gamma_+$. The spaces $M_+$ and $M$  can be considered
 as natural subsets of the mapping cylinder $\Gamma_+$,  we  will refer to $M_+$ and $M$ as the bottom 
 and the top of $\Gamma_+$  respectively. 
  
Let  $\Gamma=\Gamma(p^k)=\Gamma_+ /\z_{p^k}$,  $\gamma_+ : \Gamma_+ \lo \Gamma$  the projection
  and 
   $\Gamma_{\partial M}$ and $\Gamma_M$  the images  in $\Gamma$
   of the bottom set $M_+$  and the top set $M$  of $\Gamma_+$  under  the map $\gamma_+$.
   Note that $\alpha_+: M_+ \lo M$ induces the corresponding map $\alpha : \Gamma_{\partial M} \lo \Gamma_M$
   for which  $\Gamma$ is the mapping cylinder of $\alpha$.
   Also note that $\Gamma_M = M/\z_{p^k}$ and  $\Gamma_{\partial M}$   is   the space obtained from $M$
   by collapsing to singletons  the orbits of the action of $\z_{p^k}$ on $M$  lying in $\partial M $.
   Thus   there is a natural projection $\mu : M \lo \Gamma_{\partial M}$ and  $ \Delta^{(m+1)} \subset M$
   can be considered as a subset of $\Gamma_{\partial M}$ as well (we can identify $M$ with one of the spaces 
   $M_a\subset M_+$ and  regard $\mu $ as the projection $M_+ \lo  \Gamma_{\partial M} =M_+/\z_{p^k}$
   restricted  to $M_a$).  
   
     Let us define  a map 
 $\delta: \Gamma \lo \Delta$ by sending  $\Gamma_{\partial M}$ to $\partial \Delta$  and
 $\Gamma_M$ to the barycenter of $\Delta$ such that $\delta$ is the identity map on $\Delta^{(m+1)}$,
$\delta(x) \in \Delta'$ for every $(m+2)$-simplex $\Delta'$ of $\Delta$ and $x \in \mu(M \cap \Delta')$
and $\delta$ is linearly extended along the intervals of the mapping cylinder $\Gamma$ (in our descriptions 
we abuse notations  regarding simplexes of $\Delta$ also as subsets  of other spaces involved). 
 Denote $\delta_+ =\gamma_+ \circ \delta : \Gamma_+ \lo \Delta$.

 In the notation below
   an element  $g \in \z_{p^k}$ that   appears
 in a composition with other maps is  regarded  as a homeomorphism of the space on which $\z_{p^k}$ acts.
 
 \begin{proposition}
\label{prop2}
For every $(m+2)$-simplex $\Delta_i $ of $\Delta$  there is a map 
$r^i_+ :\delta_+^{-1}(\Delta_i) \lo \partial \Delta_i$ such that

(i)
 $r^i_+$ coincides with $\delta_+$  on  
 $ \delta_+^{-1}(\partial \Delta_i))$;

  (ii) for every triangulation of $\Gamma_+$ and  every  $g_i \in \z_{p^k}$    there is a map
   $r_+  : \Gamma_+ \lo \partial \Delta$  such that  $r_+ (\Gamma_+^{(m+2)}) \subset \Delta^{(m+1)}$ and 
  $r_+$ restricted to $ \delta_+^{-1}(\Delta_i)$ coincides with $r^i_+ \circ g_i$
  for each $i$.

\end{proposition}
{\bf Proof.} 
Recall that $\Gamma_+$ is the mapping cylinder of $\alpha_+ : M_+ \lo M$ and let
$\beta_+ : \Gamma_+ \lo M$ be the projection  to the top set $M$ of $\Gamma_+$.
Also recall that $M_+$ is the union of $p$ copies $M_a$, $a \in \z_{p^k}$, of $M$  intersecting each other at $\partial M$.

Denote $\Delta^{(m+1)}_a =\delta_+^{-1}(\Delta^{(m+1)}) \cap M_a$,  
$(\Delta_i\cap M)_a =\delta_+^{-1}(\Delta_i)  \cap M_a$,
$(\partial \Delta_i)_a=\delta_+^{-1}(\partial \Delta_i)\cap M_a$
  and $(S^m\cap M)_a =\alpha_+^{-1}(S^m\cap M) \cap M_a$.
For each $a$ and $i$ we will fix  an interval $(S_i)_a $ in 
   $(S^m\cap M)_a\cap (\Delta_i \cap M)_a$ that  connects  the sets
$\partial B_i$ and $(\partial \Delta_i)_a$.

By   (i) of Proposition \ref{prop1} the map $\beta_+$ 
 can be isotoped  relative to $\beta_+^{-1}((S^m \cap M)\cup \partial M)$ into a map 
 $\omega_+  : \Gamma_+ \lo M$  such that  $\omega_+$  and $\delta_+$ restricted
 to each $\Delta^{(m+1)}_a \subset M_+$ coincide.
 Consider the map $r$ from (ii) of Proposition \ref{prop1} and note that 
 $(r \circ \omega_+)(\Gamma_+) \subset \Delta^{(m+1)}$.
 Define the map $r_+^i$ as the map $r \circ \omega_+$ restricted   to $\delta_+^{-1}(\Delta_i)$.
\\
\\
(i) follows from (ii) of Proposition \ref{prop1}.
\\
\\
(ii)   
 For a subset $A \subset  M_+$ by the mapping cylinder 
of  $\alpha_+$ over $A$ we mean the mapping cylinder of $\alpha_+ : A \lo \alpha_+(A)$ which is
 a subset of $\Gamma_+$.
Let  $S \subset M_+$  be  the union of $(S_i)_a$ for all $i$ and $a$, and
  $B\subset M_+$ the union  of all the balls   $B_i \subset M_+$
Consider a CW-structure of $\Gamma_+$  for which the interiors of the $(m+3)$-cells are 
the interiors in $\Gamma_+$  of the mapping cylinders of
$\alpha_+$  over $(\Delta_i \cap M)_a\setminus (S_i)_a$ without the points belonging to $M_+$ and $M$. 
Thus the  $(m+2)$-skeleton of this CW-structure is the union of $M_+\cup M$ and  the mapping
cylinder of $\alpha_+$ over $\delta^{-1}(\Delta^{(2)})\cup S\cup B$.

Let be the map $\phi : M_+ \lo \partial \Delta$ be defined by the maps $r^i_+ \circ g_i$  on
each $\delta_+^{-1}(\Delta_i)$.
 Note that  the maps $\phi$   and $r \circ  \omega_+$  coincide on  
 $\delta^{-1}_+(\Delta^{(m+1)})\cup S$.
 Also note $\phi$ and  $r \circ \omega_+$ are homotopic on $B$
 by a homotopy relative to $B\cap S$. 
 Thus $r \circ \omega_+$ restricted to the union of $M \subset  \Gamma_+$ with 
 the mapping cylinder of $\alpha_+$ over
 $\delta_+^{-1}(\Delta^{(m+1)}) \cup B \cup S$
   can be homotoped
 to a map $\Phi$ such that $\phi$ and $\Phi$ restricted to  $\delta_+^{-1}(\Delta^{(m+1)}) \cup B \cup S$ coincide.
 Thus we can extend $\phi $ to a map $\phi_+$
  from  the $(m+2)$-skeleton of the CW-structure of $\Gamma_+$ to $\Delta^{(m+1)}$. 

  Then,
since $\Delta^{(m+1)}$   is contractible inside  $\partial \Delta$, we may extend $\phi_+$ to a map 
$r_+ : \Gamma_+ \lo \partial \Delta$. 
 For any triangulation  of $\Gamma_+$, the $(m+2)$-skeleton
of the  triangulation  can be pushed off the interiors of the $(m+3)$-cells of $\Gamma_+$ relative
to the $(m+2)$-skeleton of the CW-structure  of $\Gamma_+$
and  (ii) follows.
$\black$
\\
\\
Denote by $\Gamma_* $ the space obtained from $\Gamma_+$ by collapsing the fibers of $\gamma_+$
to singletons over the set $\delta^{-1}(\Delta^{(m+1)})$  and let the maps 
$\gamma_* : \Gamma_* \lo \Gamma$, $\delta_* : \Gamma_* \lo \Delta$, 
 $r^i_* : \delta_*^{-1}(\Delta_i) \lo \partial \Delta_i$ and 
 $r_* : \Gamma_* \lo \partial \Delta$ be induced by $\gamma_+$, $\delta_+$, $r^i_+$  and $r_+$ respectively and
 consider $\Gamma_*$ with the action of $\z_{p^k}$ induced by the action of $\z_{p^k}$  on
 $\Gamma_+$.  
 \begin{proposition}
 \label{prop*}
 ${}$

 (i) The conclusions of 
 Proposition \ref{prop2} hold with the subscript ``$+$''  being replaced everywhere by
 the subscript  ``$*$''.
 
 (ii)  The fixed point set of the action of $\z_{p^k}$ on $\Gamma_*$ is $(m+1)$-dimensional,
 the action of $\z_{p^k}$ is free outside the fixed point set and
  there are  triangulations of $\Gamma_*$ and $\Gamma$ and a subdivision 
 of $\Delta$ for which the action of $\z_{p^k}$ on $\Gamma_*$ is simplicial and the maps $\gamma_*$
  and $\delta$ are simplicial.
 \end{proposition}
 {\bf Proof.}
 (i) is obvious and (ii)
 can be derived from the  construction of the spaces and the maps involved. 
 $\black$
\\  \\
 Let $L$ be a finite  $(m+3)$-dimensional  simplicial complex and $\Lambda  : L \lo \Delta$ a simplicial map
 such that $\Lambda$ is $1$-to-$1$ on each simplex of  $L$. Such a map $\Lambda$ will be called 
  a {\bf sample map}. Denote by 
 $L'$ the pull-back space of the maps $\lambda$ and $\delta : \Gamma=\Gamma(p^k)  \lo \Delta$ and  by
 $\Omega: L' \lo L$ the pull-back of $\delta$.
 
 \begin{proposition}
 \label{prop3}
 The map $\Omega: L' \lo L$ induces an isomorphism  of $H_{m+3}(L'; \z_p)$ and
 $H_{m+3}(L; \z_p)$.
 \end{proposition}
 {\bf Proof.}  Recall that for every $(m+2)$-simplex $\Delta'$ of $\Delta$, $\delta^{-1}(\Delta')$
 is the mapping cylinder a  map of degree $p^k$ from $\partial \Delta'  $ to a $(m+1)$-sphere $S^{m+1}$ (the boundary of one of 
 the $(m+2)$-balls $B_i$).
 
 Also recall that  $\Gamma=\delta^{-1}(\Delta)$ is the mapping cylinder of the map $\alpha$ and from the definition of
 $\alpha$ one can also observe that 
 $H_{m+3}(\delta^{-1}(\Delta), \delta^{-1}(\partial \Delta); \z_p)=\z_p$
 and $\delta$ induces an isomorphism between 
 $H_{m+3}(\delta^{-1}(\Delta), \delta^{-1}(\partial \Delta); \z_p)$ and 
 $H_{m+3}(\Delta, \partial \Delta; \z_p)=\z_p$.
 
 And finally recall that $\delta$ is $1$-to-$1$ over the $(m+1)$-skeleton of $\Delta$. 
 Consider the long exact sequences of the pairs $(L', L'^{(m+2)})$  and $(L, L^{(m+2)})$ for the homology 
 with coefficients in $\z_p$. The facts above imply
 that $\Omega$ induces isomorphisms
 $H_{m+2}(L'^{(m+2)}; \z_p)\lo H_{m+2}(L^{(m+2)}; \z_p)$ and 
 $H_{m+3}(L', L'^{(m+2)}; \z_p)\lo H_{m+3}(L, L^{(m+2)}; \z_p)$. Then, by the $5$-lemma, $\Omega$
 induces an isomorphism $H_{m+3}(L'; \z_p)\lo H_{m+3}(L; \z_p)$ as well. $\black$
 \\
 
\begin{proposition}
\label{prop4}

 Let $G=\z_{p^k}$ act simplicially  on
a finite $(m+3)$-dimensional simplicial complex $K$ such that
  the action of $G$ is free on
$K \setminus K^{(m+1)}$. 
 Fix a  triangulation of the space $L=K/G$   for which the projection $K \lo L$ is simplicial 
 with respect to some barycentric  subdivision of the triangulation of $K$  and
 $L$ admits a sample map $\Lambda : L \lo \Delta$ to an $(m+3)$-simplex $\Delta$.
 
Then  for every $k'>k$
 there is a finite $(m+3)$-dimensional simplicial complex $K'$, a simplicial action
of $G'=\z_{p^{k'}}$ on $K'$ and a map $\omega : K' \lo K$ such that the action of $G'$ is free on
$K'\setminus K'^{(m+1)}$ and

(i)  the actions of $G$ and $G'$  agree with $\omega$ and an  epimorphism  $h : G' \lo G$.
By this we mean that $\omega(g'x))=h(g')\omega(x)$ for for every $x \in K'$ and $g' \in G'$;

(ii) there is a map $\kappa  : K' \lo K^{(m+2)}$ such that $\kappa(K'^{(m+2)}) \subset K^{(m+1)}$ and
$\kappa(\omega^{-1}(\Delta_K))\subset \Delta_K$
 for every simplex $\Delta_K$ of $K$;
 
 (iii)
  the space $L'=K'/G'$ and  the map $ \Omega: L'\lo L=L$ determined by $\omega$  can be obtained
 as the pull-back space and the pull-back map  of the sample map $\Lambda : L \lo \Delta$ and the map
 $\delta : \Gamma =\Gamma(p^{k'-k})\lo \Delta$.
 Thus, by Proposition \ref{prop3}, the map $\Omega$   induces an isomorphism
 $H_{m+3}(L'; \z_p) \lo H_{m+3}(L; \z_p)$.  
 
 (iv) Moreover,   the action of $G'$ on $K'$ is free  over
 $\Lambda'^{-1}(\Gamma\setminus \gamma_*( S))\subset L'$  
 where $\Lambda' : L' \lo \Gamma$ is the pull-back of $\Lambda$ from (iii)
  and $S \subset \Gamma_* =\Gamma_*(p^{k'-k})$
 is the fixed point set of the action of $\z_{p^{k'-k}}$ on $\Gamma_*$.
 
 \end{proposition}
 {\bf Proof.} Replacing the triangulation of $K$ by its subdivision we  assume
 that  the projection  $\pi : K\lo L$ is a simplicial map. 
 For every $(m+3)$-simplex  $\Delta_L$ of $L$  fix an $(m+3)$-simplex $\Delta_K$ of $K$ such that 
 $\pi(\Delta_K)=\Delta_L$ and denote by $K_-$ the union of the $(m+2)$-skeleton  $K^{(m+2)}$ of $K$
  with all the $(m+3)$-simplexes of $K$ that we fixed.
  Let $K'_-$ be  the pull-back space of the maps 
 $\Lambda \circ \pi|K_- : K_- \lo \Delta$ and  
 $\delta_*= \delta \circ \gamma_*: \Gamma_* =\Gamma_*(p^{k'-k})\lo \Delta$, 
 $\omega_- : K'_- \lo K$ the pull-back map of $\delta_*$ and
 $\lambda_-  : K'_- \lo \Gamma_*$ the pull back of $\Lambda\circ \pi|K_-$.

Let  $g$ be a generator  of $G'$ and  $l : G'\lo \z_{p^{k'-k}}$ an epimorphism.
 We will first  define the action of $G'$ on $\omega_-^{-1}(K^{(m+2)})$. 
 For each $(m+2)$-simplex $\Delta_L$ of $L $ define the action of $G'$ on $\omega_-^{-1}(\pi^{-1}(\Delta_L))$ as follows.
 Fix a   $(m+2)$-simplex $\Delta_K$ of $\pi^{-1}(\Delta_L)$ 
  and let  $x \in \omega_-^{-1}(\Delta_K)$. 
 Define $y=g^ t x$ for $ 1 \leq t \leq p^k-1$
 as the point $y \in K'$ such that $\omega_-(y) \in h(g^t)(\Delta_K)$ and 
 $\lambda_-(y)=\lambda_-(x)$, and for $t=p^k$ define $y=g^t x$ as
 the point $y \in \Delta_K$ such that $\lambda_-(y)=l(g)\lambda_-(x)$. 
 We do this  independently  for every $(m+2)$-simplex $\Delta_L$  of $L$ and this
 way define  the action of $g$ on $\omega_-^{-1}(K^{(m+2)})$. It is easy to see that  the action of $g$ is well-defined, 
 $g^t$ for $t={p^{k'}}$  is the identity map of  $\omega_-^{-1}(K^{(m+2)})$
  and hence the action of $g$ defines the action  of $G'$ on  $\omega_-^{-1}(K^{(m+2)})$.
  Note that

(*)    for  $g^t \in G', t=p^k $ and   $g_*=l(g)  \in \z_{p^{k'-k}}$ 
we have  that 
 $\lambda_- \circ g^t|\omega_-^{-1}(\partial \Delta_K) =g_* \circ \lambda_-|\omega_-^{-1}(\partial \Delta_K)$ 
  for every $4$-simplex $\Delta_K$ of $K$.
 
  Now we  will enlarge $K'_-$ to  a space $K'$   and  extend the action of $G'$ over 
   $K'$.
  Let $\Delta_L$ be a $(m+3)$-simplex of $L$. Recall that we fixed
   a $(m+3)$-simplex $\Delta_K$ in  $\pi^{-1}(\Delta_L)$.  For every 
  $g'=g^t \in G', 1\leq p^k-1 $ attach  to $g'(\omega_-^{-1}(\partial \Delta_K))$  a copy of 
   the space $\omega_-^{-1}(\Delta_K)$
  (which is is in its turn a copy of   $\Gamma_*$) by identifying $g'(\omega_-^{-1}(\partial \Delta_K)$   with
  $\omega_-^{-1}(\partial \Delta_K)$  according to $g'$ and for $x \in \omega_-^{-1}(\Delta_K)$ define
  $g'x$ as the as the point corresponding to $x$ in the attached space. 
We will define the action of $g^t$, $t=p^k$, on $\omega_-^{-1}(\Delta_K)$ by 
$y=g^t x, x\in \omega_-^{-1}(\Delta_K)$ such that $y \in \omega_-^{-1}(\Delta_K)$  and $\lambda_- (y)=l(g)\lambda_-(x)$.
By (*) the action of $g$ on $\omega_-^{-1}(\Delta_K)$   agrees with  the action of $g$ on $\omega_-^{-1}(K^{(m+2)})$.
 We do the above procedure  independently for every  $(m+3)$-simplex $\Delta_L$ of $L$ and
this way we define the space $K'$ and the action of $G'$ on $K'$.  We extend $\omega_-$ and $\lambda_-$ to
the maps $\omega: K' \lo K$ and $\lambda' : K' \lo \Gamma$ by 
   $\omega(g^t x)=\omega_-(x)$  and $\lambda' (g^t x) =\gamma_*   (\lambda_-(x))$
     for $x$ in a fixed $(m+3)$-simplex $\Delta_K$ and  $1\leq t \leq p^k$.

 It is easy to verify that the action of $G'$ on $K'$ and the maps $\omega$ and $\lambda'$ are
well-defined and the conclusion (i) of the proposition holds. Moreover  $\lambda'\circ g'=\lambda'$ 
for every $g' \in G'$ and  hence $\lambda'$
 defines the corresponding map $\Lambda' : L'=K'/G' \lo \Gamma$. Then $L'$ is 
the pull-back of the maps $\Lambda : L \lo \Delta$ and $\delta  : \Gamma \lo \Delta$
with $\Lambda'$ being the pull-back of $\Lambda$ and the map $\Omega : L' \lo L$ induced by $\omega$ being
the pull-back of $\delta$.
Thus, by Proposition \ref{prop3}, the conclusion (iii) of the proposition  holds as well.

Consider any triangulation of $K'$  for which
  the preimages under $\omega$ of the simplexes of $K'$ are subcomplexes of $K'$. 
Then  the map 
$r_* : \Gamma_* \lo \Delta^{(m+2)}$ provided by  Propositions \ref{prop*}  and    \ref{prop2}
for $g_i=0 \in \z_p$ for all $i$ defines 
  the corresponding map 
  $\kappa_- : K'_- \lo K^{(m+2)}$ such that $\kappa(\omega_-^{-1}(K^{(m+2)})) \subset K^{(m+1)}$. 
  The construction above  and
    Propositions \ref{prop*}  and \ref{prop2} allow us to extend 
    the map $\kappa_-$ to a map $\kappa : K' \lo K^{(m+2)}$
   satisfying the conclusion (ii) of the proposition. Recall that  the  triangulation of $K$ 
    is a subdivision of 
the original triangulation of $K$.
    Replacing $\kappa$ by its composition 
  with the simplical approximation of the identity map of $K$ with respect to the new and original triangulations 
  of  $K$ we get that
   the conclusion  (ii) of the proposition holds.
  
The rest of the conclusions of the proposition follows from
(ii) of Proposition \ref{prop*}. $\black$

\end{section}

\begin{section}{Properties related to cohomological dimension and light maps}
A CW-complex $K$ is said to be an {\bf absolute extensor}  for a space $X$ if every map $f : X \lo K$  from a closed subset 
$A \subset X$ continuously extends over $X$.  The  extension criterion of Dranishnikov 
says  that  a simply connected CW-complex is an absolute
extensor for a finite dimensional compactum $X$ if and only if 
$\dim_{H_n(K)}  X \leq n$ for every $n>1$.
\\

Recall the construction of $\Gamma$ and  note  that  for $k' >k$ there is a  natural map
$\Gamma_{\partial M}=\Gamma_{\partial M}(p^{k}) \lo 
\Gamma'_{\partial M}= \Gamma_{\partial M}(p^{k'})$ 
induced by the identity map of $M_0$ in $M_+(p^k)$ and $M_0$ in $M_+(p^{k'})$.  
Also note that  the identity map of  $M$  induces a natural map 
$\Gamma_M=\Gamma_M(p^k) \lo \Gamma'_M=\Gamma_M (p^{k'})$. 
Then these natural maps indice the corresponding map
$\tau=\tau (p^k, p^{k'})  : \Gamma= \Gamma(p^{k}) \lo  \Gamma' =\Gamma(p^{k'})$.
Also note that from the construction of $\delta: \Gamma \lo \Delta$ it follows that
for $\delta'=\delta(p^{k'}) : \Gamma' =\Gamma(p^{k'}) \lo  \Delta  $  
 we have that $\delta=\delta'  \circ \tau  $.
 \\
 
 Let  $C=C(p^k)$ be  the infinite   telescope of the maps 
 $\tau(p^t, p^{t+1}) : \Gamma (p^t) \lo \Gamma(p^{t+1})$
for $t\geq k$ and let $\delta_C : C \lo \Delta$ be the map determined by the
maps $\delta(p^t): \Gamma(p^t) \lo \Delta$.

\begin{proposition}
\label{telescope}
For every simplex $\Delta'$ of $\Delta$ the space $\delta_C^{-1}(\Delta')$ is an absolute extensor
for every finite dimensional compactum $X$ with $\p X \leq m+1$.
\end{proposition}
{\bf Proof.} 

If $\dim \Delta' \leq m+1$  then  $(\delta(p^t))^{-1}(\Delta')$ is contractible and hence the reduced homology 
   $H_n(\delta_C^{-1} (\Delta'))$  vanishes for every $n$.

If $\dim \Delta'= m+2$ then  $(\delta(p^t))^{-1}(\Delta')$ is homotopy equivalent to an $(m+1)$-sphere $S^{m+1}$
and the map $\tau (p^t , p^{t+1})$ restricted to   
 $(\delta(p^t))^{-1}(\Delta')$ and $(\delta(p^{t+1}))^{-1}(\Delta')$ acts as a map of degree $p$ between
 $(m+1)$-spheres  and
 hence  $H_{m+1}(\delta_C^{-1} (\Delta'))=\z[\frac{1}{p}] $  and    $H_n (\delta_C^{-1} (\Delta'))=0$ if $n \neq m+1$.
 
If  $\Delta' =\Delta$:  then  $(\delta(p^t))^{-1}(\Delta)=\Gamma(p^t)$.
 Recall that $\Gamma(p^t) $ is the mapping cylinder of $\alpha : \Gamma_{\partial M}(p^t) \lo \Gamma_M (p^t)$
 and hence $\Gamma(p^t)$ is homotopy equivalent to $\Gamma_M(p^t)$. Note  that $\Gamma_M (p^t)$
 is homeomorphic to $M$ and  $M$ is homotopy equivalent to the wedge of $m+3$ spheres $S^{m+1}$
 and 
 $\tau (p^t , p^{t+1})$ acts on these spheres as a map of degree $p$. Thus
 $H_{m+1}(\Gamma_C)$ is the direct sum of 
 $m+3$ copies of $\z[\frac{1}{p}]$ and $H_n(\Gamma_C)=0$ if $n \neq m+1$.
 
 Then the proposition follows from  Dranishnikov's extension criterion. $\black$
 \\
\\
A map is said to be  {\bf light or $0$-dimensional} if its fibers are $0$-dimensional.
We will say that a space $N$ has the {\bf $n$-approximation property}  if $N$ is a metric compact ANR  and
every map $f : X \lo N$ from a compactum $X$ with $\dim X \leq n$ can be arbitrarily closely
approximated by a light map. 
\\

By an {\bf extended mapping cylinder} $E$ of a map $g : N\lo N'$ we will understand 
 the mapping cylinder of  $g$ with the product $N' \times [0,1]$ attached to the mapping cylinder
 by identifying  the top  $N'$  of the mapping cylinder  with  the set $N' \times \{ 0 \}$ of the product 
 $N' \times [0, 1]$.  We will  refer to the mapping cylinder of $g$ as the proper part of $E$ and 
   to  $N'\times [0,1]$ as
 the extension part of $E$.
 
 \begin{proposition}
 \label{light}
 ${}$
 
 (i) Every compact $n$-manifold (possibly with boundary) has the $n$-approximation property.
 
 (ii) If a compactum $N$ is  covered by the interiors of closed subsets of  $N$ having
 the $n$-approximation property then $N$ has the $n$-approximation property as well.
 
 (iii) If $N$ has the $n$-approximation property then any map from an $n$-dimensional compactum 
 $f : X \lo N $ that  is light on a closed subset $A$ of $X$  
 can be arbitrarily closely
 approximated by a light map that coincides with $f$ on $A$.
 
 (iv) If $N$ has the $n$-approximation property then $N \times [0,1]$ has the $(n+1)$-approximation
 property.

 (v) If $N$ and $N'$ have the   $n$-approximation property  and $g : N \lo N'$ is a light surjective map
 then the extended mapping cylinder $E$  of $g$  has the $(n+1)$-approximation property.
 
 \end{proposition}
{\bf Proof.}
\\
(i)  follows from  the Baire category theorem applied to  mapping spaces and
 the well-known fact that a closed $n$-ball has the $n$-approximation property.
 \\
 \\
 (ii) follows from the Baire category theorem applied to mapping spaces.
 \\
 \\
 (iii) again follows from the Baire category theorem applied to mapping spaces.
 \\
 \\
 (iv) Let $f : X \lo N\times [0,1]$ be a map from a compactum $X$ with $\dim X \leq n+1$ and
 let $f=(f_N , f_I)$ be the coordinate maps $f_N : X \lo N$ and $f_I : X \lo [0,1]$ of $f$.
  Decompose $X$ into $X=A\cup B$
 with $A$ being $(n-1)$-dimensional and $\sigma$-compact and $B$ being $0$-dimensional.
 Approximate $f_N$ by a map $f'_N : X \lo N$ such that $f'_N$ restricted to $A$ is $0$-dimensional.
 Denote by $Y$ the subset of $X$ which is the union of all the non-trivial compact connected sets
 contained in the fibers of the map $(f'_N, f_I): X \lo N \times [0,1]$. Then
 $Y$ is $\sigma$-compact and $Y \cap A$ is $0$-dimensional and hence $\dim Y \leq 1$. Approximate
 $f_I$ by a map $f'_I : X \lo [0,1]$ such that $f'_I$ restricted to $Y$ is $0$-dimensional.
 Then $f'=(f'_N, f'_I) : X \lo N \times [0,1]$ is a $0$-dimensional approximation of $f$.
 \\
 \\
(v)  Let $E=E_1 \cup E_2$ be the proper and the extension parts of $E$ respectively.
Note that for every extended mapping cylinder 
 the identy map of $E$ can be arbirarily closely approximated by a map
$\phi : E \lo E$ such that $\phi^{-1}(E_1)$ is contained in the interior  of $E_1$ and
$\phi (E_2)$ is contained in the interior of $E_2$. Moreover, if $g$ is light than $\phi$ can be assumed to be light as well.
Let $ f: X \lo E$ be a map from a compactum $X$ with $\dim X \leq n+1$. By (iv) the map $ f : X \lo E$
can be approximated by a map $f' : X \lo E$ so that  $f'$ is $0$-dimensional on  a closed neighborhood
of $f^{-1}(\phi^{-1}(E_1))$ in $X$. Then taking  $f'$ to be sufficiently close to $f$ we
can also assume that $f'$ is  $0$-dimensional 
on  $(f')^{-1}(\phi^{-1}(E_1))$. Since
$\phi$ is $0$-dimensional,  the map $\phi \circ f' $  is also $0$-dimensional
on  $(f')^{-1}(\phi^{-1}(E_1))$. Clearly $X \setminus (f')^{-1}(\phi^{-1}(E_1)) \subset (f')^{-1}(\phi^{-1}(E_2))$.
Then,  by (iii) and (iv),  the map $\phi \circ f' : X \lo E$  restricted to $ (f')^{-1}(\phi^{-1}(E_1))$ can be extended
arbitrarily closely to $\phi \circ f'$ to a light map $f'' : X \lo E$.
$\black$
\\
\\
Let us return again to  the construction of $\Gamma=\Gamma(p^k)$ and $\delta : \Gamma \lo \Delta$. 
For every $(m+2)$-simplex  $\Delta_i$ of $\Delta$ the set $\delta^{-1}(\Delta_i) \subset \Gamma_{\partial M}$ 
can be represented
as  the mapping cylinder
of a map of degree $p^k$ from $\partial \Delta'$ to the sphere 
$S^{m+1}_i=\gamma_+(\partial B_i)\subset \Gamma_{\partial M}$,
$\partial B_i \subset M_+$.
 Enlarge this mapping cylinder
to  the extended one by attaching $ \partial S^{m+1}_i \times [0,1]$ to the top $S^{m+1}_i$ of the mapping cylinder.
Let us also attach $S^{m+1}_i  \times [0,1]$ to  
$S^{m+1}_i  =\gamma_+ (\partial B_i)\subset  \Gamma_M, \partial B_i \subset M$. 
Denote by $\Gamma_{\partial M}^E$ the space obtained from $\Gamma_{\partial M}$ after
 the attachments  for all  $i$ and denote by $\Gamma_M^E$ the space obtained from 
 $\Gamma_M$ after the attachments for all $i$. 
 Note that $\Gamma^E_M$ is an $(m+2)$-manifold with boundary 
 (actually both $\Gamma^E_M$ and  $\Gamma_M$  are homeomorphic to $M$).
 Clearly the map
 $\alpha : \Gamma_{\partial M} \lo \Gamma_M$ extends to
 $\alpha_E : \Gamma_{\partial M}^E \lo \Gamma_M^E$ by the identity map on the attached parts.
 Recall that $\Gamma$ is the mapping cylinder of $\alpha$ and denote by
 $\Gamma_E$ the extended mapping cylinder $\alpha_E $. Note that the  natural retractions 
 of $\Gamma^E_{\partial M}$ to $\Gamma_{\partial M}$ and
 of $\Gamma^E_M$ to $\Gamma_M$ induce the corresponding retraction from 
 the mapping cylinder of $\alpha_E$ to $\Gamma$, and composing the last retraction with
 the natural retraction from the extended mapping cylinder $\Gamma_E $ to the mapping cylinder of $\alpha_E$
 we obtain the retraction $\gamma_E : \Gamma_E \lo \Gamma$.  Denote $\delta_E =\delta \circ \gamma_E : \Gamma_E \lo \Delta$.
 
\begin{proposition}
\label{extended-gamma}
${}$

(i) $\dim \Gamma_E =m+3$ and $\dim \gamma_E^{-1}(F)=m+1$ where
$F$ is the image under $\gamma_*$ of the fixed point set of the action of $\z_{p^k}$ on $\Gamma_*$.

(ii) $\Gamma_E $ has the $(m+3)$-approximation property. Moreover
for every $(m+2)$-simplex $\Delta'$ of $\Delta$ the space 
$\delta_E^{-1}(\Delta')  \subset \Gamma_E$ has the $(m+2)$-approximation
property.
\end{proposition}
{\bf Proof.}

(i)
$\dim \Gamma_E =m+3$ is obvious.  Note that  $F$
consists of $F_1 =\delta^{-1}(\Delta^{(m+1)})$, 
 $F_2 =$ the union of  $(S^{m+1}_i \cap S^m) \times [0,1]$ for all $i$  and 
 $F_3=S^m\subset \Gamma_M$ where $S^m$ is the fixed point set of the action of $T$ on $S^{m+2}$.
Then $\gamma_E^{-1}(F)$ can be represented as the union of the following parts 
homeomorphic  to $F^E_1=F_1$,  $F^E_2=F_2\times [0,1] $ and  
$F^E_3 = ( F_2 \cup F_3) \times [0,1]$.
 Note that $\dim F_1 =m+1$ and $\dim F_2=m$ and $\dim F_3 =m$. Then
 $\dim F^E_1=m+1$, $\dim F^E_2=m+1$ and $\dim F^E_3=m+1$ and hence $\dim \gamma_E^{-1}(F)=m+1$.
\\

(ii) By Proposition \ref{light} the spaces $\delta^{-1}_E (\Delta')$,
 $\Gamma_{\partial M}^E $ and $\Gamma_M^E$
have the $(m+3)$-approximation property and hence, again by Proposition \ref{light}, the space
$\Gamma_E$ has the $(m+3)$-approximation property. $\black$
\\

Let $L'$ be a CW-complex and $L$ a simplicial complex.
We say that  a map $\Omega : L' \lo L$ is {\bf combinatorial} 
if the preimage under $\Omega$ of every simplex of $L$ is a subcomplex of $L'$.
We also say that for a map   $f  : X \lo L$ a map $f' : X \lo L'$ is a {\bf combinatorial lifting} of $f$ if
$\Omega(f'(f^{-1}(\Delta_L))) \subset \Delta_L$ for every simplex $\Delta_L$ of $L$.
It is easy to see that if for a space $X$ and  a combinatorial map $ \Omega : L'\lo L$ we have that 
$\Omega^{-1}(\Delta_L)$ is an absolute extensor for $f^{-1}(\Delta_L)$ for every simplex $\Delta_L$ of $L$
then any map $f : X \lo L$ admits a combinatorial lifting $f' : X \lo L'$.  
Moreover, if $A$ is a closed subset of $X$ and $f'_A : A \lo L'$ is a combinatorial lifting of $f$
restricted to $A$ then $f'_A$ extends to a combinatorial lifting $f' : X \lo L'$  of $f$.
\\

Let $\Lambda : L \lo \Delta$ be a sample map of  a finite  $(m+3)$-dimensional simplicial complex $L$
to an $(m+3)$-simplex $\Delta$.

\begin{proposition}
\label{next-k}
Let $f : X \lo L$ be  a map  from a finite dimensional compactum $X$ with $\p X \leq m+1$.
Then for   every  natural number $k$ there is a natural number $k'> k$ such that 
 for  the pull-back space $L'$ of the maps $\Lambda$ and 
$\delta'=\delta(p^{k'}) : \Gamma'=\Gamma(p^{k'}) \lo \Delta$  there is a combinatorial lifting 
$f' : X \lo L'$   of $f$ with respect to the pull-back map $\Omega: L' \lo L$
of $\delta'$.
\end{proposition}
{\bf Proof.}
Consider the infinite telescope $C(p^k)$ with $\Gamma(p^k)$ being naturally embedded in $C(p^k)$ and
let $L_C$ be the pull-back space of the  sample map  $\Lambda : L \lo \Delta$ and
$\delta_C : C(p^k) \lo \Delta$ and   let $\Omega_C : L_C \lo L$ 
and $\Lambda_C : L_C \lo C(p^k)$ 
be the pull-back maps   of $\delta_C$ and $\Lambda$ respectively.

Then, by Proposition \ref{telescope}, the map $f$ admits a combinatorial lifting 
$f_C : X  \lo C(p^k)$ to $C(p^k)$ with respect to $\Omega_C$. Since $X$ is compact there is  $k'>k$
such that $\Lambda_C(f_C(X)) \subset C(p^k)$ is contained in the finite part $C'$ of the telescope 
$C(p^k)$ that ends at $\Gamma (p^{k'})$. Let $\pi : C' \lo \Gamma(p^{k'})$ be the natural
projection and let $\Pi : L_C \lo L'$ be the map induced by $\pi$. Then
$f'=\Pi \circ f_C : X \lo L'$ is  the required  combinatorial lifting of $f$. $\black$
\\

Let $L$ be a finite $(m+3)$-dimensional  simplcial complex and $\Lambda : L \lo \Delta$ a sample
map to an $(m+3)$-simplex  $\Delta$.
Consider the pull-back spaces $L'$ and $L'_E$ of the map $\Lambda$ and
the maps  $\delta:\Gamma \lo  \Delta$ and $ \delta_E=\delta \circ \gamma_E : \Gamma_E \lo \Delta$
respectively
and let
$\Omega : L' \lo L$ and $\Omega_E : L'_E \lo L$  be the pull-backs of $\delta$ and
$\delta_E$ respectively and let $r_E : L'_E \lo L'$ be the  retraction  induced by 
the retraction $\gamma_E : \Gamma_E \lo \Gamma$.

\begin{proposition}
\label{next-light}
 Assume that  a map $f : X \lo L$  from a compactum $X$ is 
such that   $\dim X \leq m+3$,  $\p X \leq m+1$  and $\dim f^{-1}(\Delta_L) \leq \dim \Delta_L$
for every simplex $\Delta_L$ of $L$. Also assume that $f$
admits a combinatorial lifting $f' : X \lo L'$ with respect to $\Omega$.  Then
$f$ also admits a combinatorial lifting $f'_E : X \lo L'_E$ with respect to $\Omega_E$ such that 
$f'_E$ is a light map.
\end{proposition}
{\bf Proof.} The proposition follows from (ii) of Proposition \ref{extended-gamma}. $\black$

\end{section}

\begin{section}{Proof of  Theorem  \ref{main-th}}

\begin{proposition}
\label{main-prop}
Let $X$ be a compactum with $\dim X \leq m+3$ and $\p X \leq m+1$. Then
there is a compactum $Y$ and an action of $A_p$ on $Y$  such that $\dim Y \leq m+3$,
$\dim Y/A_p \leq m+1$  and $X$ admits a light map to $Y/A_p$.
\end{proposition}
{\bf Proof.}
We will construct by induction finite simplicial complexes $K_i$,  a simplicial  action of $\z_{p^{k_i}}$  on $K_i$,
bonding maps $\omega_{i+1} : K_{i+1} \lo K_i$, $\kappa_{i+1}: K_{i+1} \lo K_i$
 and   maps $f_i : X \lo L_i=K_i /\z_{p^{k_i}}$
 such that  $k_{i+1} \geq  k_i$, 
the map $\omega_{i+1}$ agrees with the actions of $\z_{p^{k_{i+1}}}$ and $\z_{p^{k_i}}$ on
$K_{i+1}$ and $K_i$ respectively  with respect to an epimorphism $\z_{p^{k_{i+1}}}\lo \z_{p^{k_i}}$.

Set $L_0=K_0$  to be an  $(m+3)$-dimensional sphere, $f_0 : X \lo L_0$ to be  any light map
and let the trivial group $\z_{p^0}=0$  trivially act on $K_0$.
Assume  that the construction is completed for $i$.  To  proceed to $i+1$  we will choose one of  the following
procedures.
\\

{\bf Procedure  1.} Let $K=K_i, L=L_i, k=k_i$ and  $\Lambda: L \lo K$  the sample map 
 satisfying the assumptions of  Proposition \ref{prop4}.
 Assume that $f_i : X \lo L_i$ is light,
 apply Proposition \ref{next-k} for $f=f_i$  to produce $k'$ and set
 $k_{i+1}=k_i +k'$. 
  Apply Proposition \ref{prop4} with  $K=K_i$, $L=L_i$, $k=k_i$, $k'=k_{i+1}$ 
   to construct $K_{i+1}=K'$, $L_{i+1}=L'$, $\omega_{i+1}=\omega$, 
   $\kappa_{i+1}=\kappa$  and
  the action of $\z_{p^{k_{i+1}}}$ on $K_{i+1}$ satisfying the conclusion of Proposition \ref{prop4}.
  By   (iii) of Proposition \ref{prop4} and Proposition \ref{next-k} the map $f_i$ admits a combinatorial lifting
  $f_{i+1} : X \lo L_{i+1}$ with respect to $\Omega=\Omega_{i+1} : L_{i+1} \lo L_i$.
  \\
  
  {\bf Procedure 2.}  Assume that   we  proceeded from $i-1$ to $i$ according to Procedure 1. 
  Let $L=L_{i-1}$ and  $\Lambda : L=L_{i-1} \lo \Delta$ be the sample map
  used in Procedure 1 for proceeding from $i-1$ to $i$.  Since the map $f_{i-1}$  in Procedure 1 was assumed to be  light
   we have that $\dim f^{-1}_{i-1}(\Delta_L) \leq \dim \Delta_L$ for every simplex $\Delta_L$ of $L=L_{i-1}$. 
  Apply Proposition \ref{next-light}  for  $f=f_{i-1}, f'=f_i$, $L=L_{i-1}$ to construct 
  $L_{i+1}=L'_E$ and    a light map $f_{i+1}=f'_E$.
  Define $K_{i+1}$ as the pull-back space  of the projection $K_i \lo L_i = K_i/\z_{p^{k_i}}$ and
and 
  $\Omega_{i+1}=\Omega_E : L_{i+1}=L'_E \lo L_i =L$, and
  define $\omega_{i+1} : K_{i+1} \lo K_i$ to be the pull-back map of $\Omega_{i+1}$.
  Set $k_{i+1}=k_i$ and 
  define  the action of $\z_{p^{k_{i+1}}}$  
  on $K_{i+1}$
  as  the pull-back action of $\z_{p^{k_i}}$ on $K_i$.
   By 
   (i) of Proposition \ref{extended-gamma}, 
    (iii) and (iv)  of Proposition \ref{prop4} and (ii) of  Proposition \ref{prop*}
  one can find a triangulation of $L_{i+1}$ whose preimage under the projection
  $K_{i+1} \lo L_{i+1} =K_{i+1} /\z_{p^{k_{i+1}}}$  provides a triangulation
  of $K_{i+1}$ for which the action of $\z_{p^{k_{i+1}}}$  on $K_{i+1}$ is simplicial and
  free outside the $(m+1)$-skeleton of $K_{i+1}$. Set $\kappa_{i+1} : K_{i+1} \lo K_i$
  to be s simplicial approximation of $\omega_{i+1}$.
\\  
  
Carry out the construction  applying Procedure 1 for $i=2t$ and  Procedure 2 for $i=2t+1$.
Let $\Omega_{i+1} : L_{i+1} \lo L_i$ be the map determined by  
$\omega_{i+1}: K_{i+1} \lo K_i$. Denote  $Y=\invlim (K_i, \omega_i)$ and
consider $Y$ with the action of $A_p=\invlim \z_{p^{k_i}}$ determined by the actions
of $\z_{p^{k_i}}$ on $K_i$ for all $i$.  Then $Y/A_p =\invlim (L_i, \Omega_i)$.
Note that    we can replace the triangulation of $K_{2t}$
by any of its barycentric subdivisions and assume that the triangulation of $K_{2t}$  is as fine as we wish.
Then (ii) of Proposition \ref{prop4} guarantees  that $\dim Y \leq m+1$ and
Procedure 2  allow us to acheive that there is a map $f : X \lo Y/A_p$ determined by  the maps $f_i$
such that  for each $t>0$ the map $f$ followed by the projection of $Y/A_p$ to
$L_{2t}$ is as close to $f_{2t}$ as  we wish. Recall that $f_{2t}$ is light. Then the construction
can be carried out so that $f$ is light map. Clearly $\dim Y/A_p \leq m+3$
and the proposition follows. $\black$
\\
\\
{\bf Proof of Theorem \ref{main-th}.}  Let $f : X \lo Y/A_p$  be the light map produced 
by Proposition \ref{main-prop} for $m=n-1$. 
Denote by  $Z$ be the pull-back space of  $f$ and the projection $Y\lo Y/A_p$ and 
endow $Z$  with the pull-back action of $A_p$ on $Y$. Then $X=Z/A_p$ and 
the pull-back map $Z \lo Y$ of $f : X \lo Y/A_p$ is light since $f$ is light. Thus we have $\dim Z \leq \dim Y \leq m+1=n$.
$\black$

\end{section}
\begin{section}{Concluding remarks}
 Yang \cite{yang} showed that for an action of $A_p$ on a compactum $Z$ we have
 $\dim_{\z} Z/A_p \leq \dim_{\z} Z +3$, $\dim_{\z_p} Z/A_p \leq \dim_{\z_p} Z +2$ and
  $\dim_{\z_{p^\infty}} Z/A_p \leq \dim_{\z_{p^\infty}} Z +2$.
  It is still an open problem if the equality $\dim_{\z} Z/A_p = \dim_{\z} Z +3$ can be achieved.
 One can show that both  Raymond-Williams' example and  the example in \cite{levin-on-rw}
 satisfy $\dim_{\z_p} Z/A_p =\dim_{\z_p} Z +2$. Theorem \ref{main-th} shows that
 the gap $2$ can be also achieved for $\dim_{\z_{p^\infty}}$.
 
 \begin{corollary}
 There is an action of $A_p$ on an $n$-dimensional compactum $Z$, $n\geq 2$, such that 
 $\dim_{\z_{p^\infty}} Z/A_p=\dim_{\z_{p^\infty}} Z +2=n+1$.
 \end{corollary}
 {\bf Proof.} Take any compactum $X$ with $\dim X=n+2$, $\p X =n-1$ and $\dim_{\z_{p^\infty}} X =n+1$.
 Note that the dimensional restrictions on $X$ satisfy Bockstein's inequalities and hence, by
 Dranishnikov's realization theorem, such a compactum $X$ exists. Apply Theorem \ref{main-th}
 to construct an $n$-dimensional compactum $Z$ and an action of $A_p$ on $Z$ such that
 $X=Z/A_p$.  Recall that $\p Z =\p X =n-1$. Then, by Bockstein inequalities, we have that
  $\dim_{\z_{p^\infty}} Z = n-1$. $\black$
 \\
 
  A similar argument based on Yang's relations and Bockstein's inequalities also shows that
 
 \begin{corollary}
 No  $(n+3)$-dimensional compactum  $X$ with $\p X \leq n-1$ can be obtain as the orbit space
 of an action of $A_p$ on an $n$-dimensional compactum $Z$.
 \end{corollary} 
 
 We will end with a couple of related problems.
 
 \begin{problem}
  Does there exist a  $3$-dimensional compactum that  can be obtained as the orbit space
 of an action of $A_p$ on a $1$-dimensional compactum? 
 Can any $3$-dimensional compactum $X$ with $\p  X =1$ be obtained as the orbit space
 of an action of $A_p$ on a $1$-dimensional compactum? 

 \end{problem}
 
 \begin{problem} Let $n \geq 1$.
 Does there exist an  $(n+2)$-dimensional compactum $X$ that can be obtained as the orbit space
 of an action of $A_p$ on an $n$-dimensional compactum such that the action of $A_p$ is  free over
 a subset  $X'$ of  $X$ of $\dim X \setminus X' \leq n-1$?
 Can any $(n+2)$-dimensional compactum $X$ with $\p  X \leq n$  be obtained as the orbit space
 of an action of $A_p$ on an $n$-dimensional compactum such that the action of $A_p$ is  free over
 a subset of $X'$ of  $X$ with $\dim X \setminus X' \leq n-1$?
 
 \end{problem}
  Note that the action of $A_p$ constructed in Theorem \ref{main-th} is not free over  a subset of $X$ of $\dim=n$.

\end{section}

Michael Levin\\
Department of Mathematics\\
Ben Gurion University of the Negev\\
P.O.B. 653\\
Be'er Sheva 84105, ISRAEL  \\
 mlevine@math.bgu.ac.il\\\\
\end{document}